\tolerance=10000
\magnification=\magstephalf
\def\nagy{\font\caps=cmcsc10\caps}
\def\kis{\font\kisf=cmr5\kisf}

\def\sethoff{\hoffset=0.5 true in}
\sethoff  
\hsize=6.0 true in
\voffset=.3 true in
\vsize=8.5 true in  

\newif\ifduplexpr
\duplexprfalse

\def\duplex{\hoffset=0.25 true in \duplexprtrue \marginfalse}


\def\uj {\bigskip \rm}


\def\redef#1#2{\expandafter\ifx\csname #1\endcsname\relax
    \expandafter\edef\csname #1\endcsname{#2} 
    \else \message{redefinition of '\string#1'
    } \fi }

\newif\ifexist

\def\testfile#1{
\openin 0=#1
\ifeof 0
\existfalse
\else
\existtrue
\fi
\closein 0
}


\testfile{cite.inc}
\ifexist
\input cite.inc
\else
\message{!!!!!!!!!! One more pass needed for references !!!!!!!!!}
\fi


\immediate\openout 0=cite.inc


\immediate\openout 2=conten.inc
\def\writeitem#1#2{\write2
  {\string
  \line{#1\string\nagy{}#2\string\rm\string\leaderfill\string\quad\folio}}}


\def\boxit#1{\vbox{\hrule\hbox{\vrule#1\vrule}\hrule}}


\newif\ifmargin
\margintrue


\newcount\sorszam \sorszam=0
\def\sorszaminc {\advance\sorszam by1 
\ifnum \section=0 \else \the\section.\fi
\the\sorszam. }

\def\authorstr{}
\def\shorttitlestr{}

\def\author#1{\def\authorstr{#1} 
 \bigskip
 \centerline{#1}
}

\def\abstract#1{\bigskip
 {\advance\leftskip by 1in
 \centerline{\nagy Abstract}\medskip\noindent
 {\narrower#1}}
}

\def\shorttitle#1{\def\shorttitlestr{#1}}

\def\oldalszam{
\headline={
 \nagy
 \ifnum \folio > 1
 \ifodd \folio {\hfil \shorttitlestr \hfil\folio} 
 \else
 {\folio\hfil \authorstr \hfil}
 \fi
 \else \hfil
 \fi
}
}

\nopagenumbers
\oldalszam


\newcount\section \section=0
\def\newsection#1{
      \goodbreak
      \advance\section by1 \sorszam=0 
      \medskip\bigskip\centerline{\nagy \the\section. #1}
      \nobreak\nobreak
     \writeitem{\the\section.  }{#1}
     }
\newbox\labelboxx
\def\ugor{\goodbreak
\bigskip
\ifmargin
\hskip-1in {\box\labelboxx}
\vskip-\baselineskip
\fi
\noindent\bf}
\def\skippy{\enskip\enskip}

\def\definition { \ugor Definition \sorszaminc \rm \ }   
\def\theorem {    \ugor Theorem \sorszaminc \sl \ }   

\def\lemma {      \ugor Lemma \sorszaminc \sl \ }

\def\proof { \medskip {\noindent \it Proof.\skippy}\rm}
\def\example {    \ugor Example \sorszaminc \rm \ }

\def\nullbox{\setbox0=\null \ht0=5pt \wd0=5pt \dp0=0pt \box0}
\def\eop {\hfill \boxit{\nullbox}\goodbreak}


\def\label#1{  
\ifmargin
\vskip0in\hskip-1in {\kis#1}
\vskip-\baselineskip
\fi
\immediate
\write0{\string\redef{\string#1}{{\string\rm\the\section.\the\sorszam}}}}

\def\mark#1{  
\advance\sorszam by 1
\setbox\labelboxx=\hbox{\kis #1}
\immediate
\write0{\string\redef{\string#1}{{\string\rm
\ifnum \section=0 \else\the\section.\fi
\the\sorszam}}}
\advance\sorszam by -1
}

\def\labelref#1{ 
\immediate
\write0{\string\redef{\string#1}{{\string\rm\the\sorszam}}}}

\def\cite#1{\expandafter\ifx\csname#1\endcsname\relax
    ???
    \message{!!!!!!! Missing reference '\string#1' !!!!!!!}
    \else
\csname#1\endcsname
\fi}


\def\title#1{\shorttitle{#1}\centerline{\nagy #1}}


\def\references{
      \goodbreak
      \medskip\bigskip\centerline{\nagy References}
      \nobreak\nobreak
     \writeitem{\ }{References}
}

\def\i#1#2#3#4{\itemitem{\hbox to .5in{#1\hfil}}{#2, {\it #3}, #4}}



\def\smallmatrix#1{\null\,\vcenter{\baselineskip=0.7\baselineskip
    \ialign{\hfil$\scriptstyle##$\hfil&&$\,\,$\hfil$\scriptstyle##$\hfil\crcr
      \mathstrut\crcr\noalign{\kern-\baselineskip}
      #1\crcr\mathstrut\crcr\noalign{\kern-\baselineskip}}}\,
      \normalbaselines}

\input xy.tex
\xyoption{matrix}
\xyoption{curve}
\xyoption{arrow}
\xyoption{ps}
\xyoption{cmtip}

\duplex

\title{A NOTE ON PRIMITIVE EQUIVALENCE}

\author{N\'andor Sieben}
\shorttitle{Note on Primitive Equivalence}

\footnote{}{1991 {\it Mathematics Subject Classification.} Primary 46L05, 
Secondary 05C20 }

\abstract{Primitive equivalence of graphs and matrices was used by 
Enomoto, Fujii and Watatani to classify Cuntz-Krieger algebras of 
$3\times 3$ irreducible matrices.  In this paper it is shown that the 
definition of primitive equivalence can be simplified using primitive 
transfers of matrices that involve only two rows of the matrix.  } 

\newsection{Introduction}

\uj Primitive equivalence of graphs and matrices was used by 
Enomoto, Fujii and Watatani [EFW] to classify Cuntz-Krieger algebras 
of $3\times 3$ irreducible matrices.  It was shown by Drinen and the author 
[DS] that a graph and its primitive transfer have isomorphic 
groupoids and therefore isomorphic $C^{*}$-algebras.  Franks [Fra, 
Corollary 2.2] used a similar operation to find a canonical form for 
the flow equivalence class of a matrix.  His definition is more 
general but involves only two columns of a matrix.  The purpose of 
this paper is to show that the definition of primitive equivalence can 
be simplified using primitive transfers that involve only two rows of 
the matrix.  This fact can be used to simplify proofs about primitive 
equivalence in [EFW] and in [DS]. The author thanks Doug Drinen and 
John Quigg for their help.

\newsection{Preliminaries}

\uj A {\it digraph} $E$ is a pair $(E^0,E^1)$ of possibly infinite sets 
where $E^1\subset E^0\times E^0$.  $E^0$ is 
called the set of vertices and $E^1$ is called the set of edges.  We say 
that the edge $e=(v,w)$ starts at vertex $v$ and ends at vertex $
w$.  If 
$v=w$ then $e$ is called a loop. Note that a digraph has no multiple edges
but it can have loops.
The {\it vertex matrix} $A$ of a graph $
E$ is 
an $E^0$ by $E^0$ matrix such that 
$$A_{v,w}=\cases{1&if $(v,w)\in E^1$\cr
0&else.\cr}
$$
There is a bijective correspondence between digraphs and 
$0-1$ square matrices. A {\it connected component\/} of a digraph 
is a maximal subgraph in which every two vertices can be connected 
by an undirected path. 

If $A$ is a matrix then we denote the $i$-th row of $
A$ by $A_i$ 
and we use the notation $E_j=(0,\dots,0,1,0,\dots,0)$ for a row which 
has a $1$ in the $j$-th column and $0$ in all the other columns.  

Let $A$ be the vertex matrix of a digraph $E$, 
that is, a $0-1$ square matrix.  If 
$$A_p=A_{m_1}+\cdots +A_{m_s}+E_{k_1}+\cdots +E_{k_r}$$
for some distinct $k_1,\dots,k_r,m_1,\dots,m_s$ and $p\not \in \{
m_1,\dots,m_s\}$ then the 
$0-1$ matrix $B$ defined by 
$$B_{i,j}=\cases{A_{i,j}&if $i\neq p$\cr
1&if $i=p$ and $j\in \{m_1,\dots,m_s,k_1,\dots,k_r\}$\cr
0&else\cr}
$$
is called by [EFW] a {\it primitive transfer\/} of $A$ at $p$ (see also [DS]). 

We call the number of elements in $M=\{m_1,\dots,m_s\}$ the {\it size\/} of the 
primitive transfer.  Alternatively, $B$ can be defined as 
$$B_i=\cases{A_i&if $i\not =p$\cr
A_p-\sum_{m\in M}A_m+\sum_{m\in M}E_m&if $i=p${\rm .}\cr}
$$
The digraph $F$, whose vertex matrix is $B$, is also called a 
{\it primitive transfer\/} of $E$.  

\mark{ptex}

\example If 
$$A=\left(\smallmatrix {1&1&0&1&1&1&1&1\cr
0&0&0&0&0&0&0&0\cr
0&1&0&0&0&0&0&0\cr
0&0&0&1&1&0&0&0\cr
0&0&0&0&0&0&0&0\cr
1&0&0&0&0&0&0&0\cr
0&0&0&0&0&1&0&1\cr
0&0&0&0&0&0&1&0\cr}\right)\quad\hbox{\rm and}\quad\quad B=\left(\smallmatrix {0&0&1&1&0&1&1&1\cr
0&0&0&0&0&0&0&0\cr
0&1&0&0&0&0&0&0\cr
0&0&0&1&1&0&0&0\cr
0&0&0&0&0&0&0&0\cr
1&0&0&0&0&0&0&0\cr
0&0&0&0&0&1&0&1\cr
0&0&0&0&0&0&1&0\cr}\right)$$
then $B$ is a primitive transfer of $A$ since $A_1=A_3+A_4+A_6+A_
7+A_8$ 
and 
$$B_1=A_1-A_3-A_4-A_6-A_7-A_8+E_1+E_3+E_4+E_6+E_7+E_8.$$
The corresponding graphs are:
$$
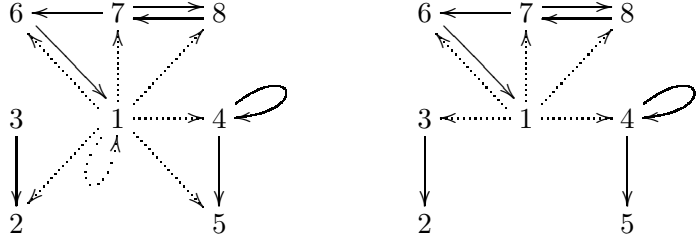
\xymatrix{
6\ar@<2pt>[dr]          &7\ar[l]
			  \ar@<2pt>[r]          &8\ar@<2pt>[l]  \cr
3\ar[d]                 &1\ar@{.>}@<-2pt>[ld]
			  \ar@{.>}@<2pt>[ul]
			  \ar@{.>}[u]
			  \ar@{.>}[ur]
			  \ar@{.>}[r]
			  \ar@{.>}[rd]
			  \ar@{.>}@(ld,d)[]
						&4\ar[d] 
						  \ar@(ru,r)[]  \cr
2		      	&                       &5              \cr
}
\qquad\qquad\qquad
\xymatrix{
6\ar@<2pt>[dr]          &7\ar[l]
			  \ar@<2pt>[r]          &8\ar@<2pt>[l]  \cr
3\ar[d]                 &1\ar@{.>}[l]
			  \ar@{.>}@<2pt>[ul]
			  \ar@{.>}[u]
			  \ar@{.>}[ur]
			  \ar@{.>}[r]
						&4\ar[d]       
						  \ar@(ru,r)[]  \cr
2                 	&                       &5              \cr
}
$$ 

\definition If $A$ and $B$ are the vertex matrices of the digraphs 
$E$ and $F$ then $A$ and $B$ are called {\it primitively equivalent\/} if there exist 
matrices $A=C_1,\dots,C_i,\dots,C_n=B$ such that for all $1\le i\le 
n-1$ one 
of the following holds:  
\item{(i)}{$C_i$ is a primitive transfer of $C_{i+1}$;} 
\item{(ii)}{$C_{i+1}$ is a primitive transfer of $C_i$;} 
\item{(iii)}{$C_i=PC_{i+1}P^{-1}$ for some permutation matrix $P$.} 

\noindent
The digraphs $E$ and $F$ are also called {\it primitively equivalent}.

\newsection{Main theorem}

\uj The purpose of this paper is to show that in the definition of 
primitive equivalence we could consider only size $1$ primitive 
transfers.  To see this we are going to show that a primitive 
transfer can be replaced by a sequence of size $1$ primitive transfers.  
First we need a few tools.  

\mark{graphdef}

\definition Let $E$ be a digraph and let $B$ be a 
primitive transfer of the vertex matrix $A$ of $E$, corresponding to the 
equation $A_p=\sum_{m\in M}A_m+\sum_{k\in K}E_k$.  The {\it graph of the primitive }
{\it transfer\/} is the subgraph of $E$ induced by $M$.  

\uj Note that the graph of a primitive transfer has no vertex with 
more than one incoming edge since $A_p$ contains only $0$'s and $
1$'s, hence 
no two $A_m$'s can have a $1$ at the same location.  

\mark{gex}

\example The graph of the primitive transfer in Example \cite{ptex} has
three connected components:
$$
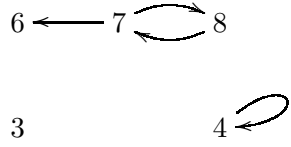
\xymatrix{
6                       &7\ar[l]
			  \ar@/^/[r]            &8\ar@/^/[l]    \cr
3			&			&4             
						  \ar@(ru,r)[]  \cr
}
$$ 

\mark{edgelemma}

\lemma If $G$ is the graph of a primitive transfer with size $s$ from the 
matrix $A$ to $B$ and there is an edge $(l,n)$ in $G$ which is not a loop, 
then there is a matrix $C$ such that $A$ is a size $1$ primitive transfer 
of $C$ and $B$ is a size $s-1$ primitive transfer of $C$.  If $H$ is the graph 
of the latter primitive transfer then $(l,m)\in H^1$ whenever $(n
,m)\in G^1$.  

\proof Let $B$ be the primitive transfer of $A$ at $p$ corresponding to the 
equation $A_p=\sum_{m\in M}A_m+\sum_{k\in K}E_k$.  Note that $M$ and $
K$ are disjoint 
and $p\notin M$.  Define 
$$C_i=\cases{A_i&if $i\neq l$\cr
A_l+A_n-E_n&if $i=l$.\cr}
$$
$C$ is a $0-1$ matrix since $A_l+A_n$ is a $0-1$ row and $
A_{l,n}=1$. 
If $J=\{j:A_{l,j}=1\hbox{\rm \ and $j\neq n$}\}$ then the equation 
$$C_l=A_n+\sum_{j\in J}E_j$$
determines a primitive transfer of $C$ at $l$. This primitive transfer is 
$A$ since
$$A_i=\cases{C_i&if $i\neq l$\cr
C_l-C_n+E_n&if $i=l.$\cr}
$$
The equation 
$$\eqalign{C_p&=A_p=\sum_{m\in M}A_m+\sum_{k\in K}E_k\cr
&=\sum_{m\in M\setminus \{l\}}C_m+(C_l-C_n+E_n)+\sum_{k\in K}E_k\cr
&=\sum_{m\in M\setminus \{n\}}C_m+\sum_{k\in K\cup \{n\}}E_k\cr}
$$
determines another primitive transfer of $C$ at $p$.  This 
primitive transfer is $B$ since $B_i=A_i=C_i$ for $i\neq p$ and 
$$\eqalign{B_p&=A_p-\sum_{m\in M}A_m+\sum_{m\in M}E_m\cr
&=C_p-\sum_{m\in M\setminus \{n,l\}}C_m+\sum_{m\in M\setminus \{n
\}}E_m-(A_l+A_n-E_n)\cr
&=C_p-\sum_{m\in M\setminus \{n\}}C_m+\sum_{m\in M\setminus \{n\}}
E_m.\cr}
$$
The last part of the lemma follows from the construction of $C$.  
\eop

\uj Note that in the previous lemma, $H$ is connected if $G$ is.  This 
follows from the last sentence of the lemma and the fact that $(l
,n)$ 
is the only edge in $G$ that ends at $n$.  

\mark{complemma}

\lemma If $G$ is the graph of the primitive transfer from the matrix $
A$ 
to $B$ and $F$ is one of several connected components of $G$, then there 
is a matrix $D$ such that $D$ is the primitive transfer of $A$ with graph 
$F$, and $B$ is the primitive transfer of $D$ with graph $H:=G\setminus 
F$.  

\proof Let the original primitive transfer correspond to the equation 
$A_p=\sum_{m\in M}A_m+\sum_{k\in K}E_k$. Note that $M=F^0\cup H^0$ and $
F^0\cap H^0=\emptyset$.  If 
$J=\{j:A_{h,j}=1\hbox{\rm \ for some }h\in H^0\}$ then the equation 
$$A_p=\sum_{m\in F^0}A_m+\sum_{k\in K\cup J}E_k$$
determines a primitive transfer of $A$ at $p$.  Note that $F^0$ and $
K\cup J$ 
are disjoint since $M\cap K=\emptyset$ and there is no edge from $
H^0$ to $F^0$ in $G$.  
Let $D$ be this primitive transfer of $A$ at $p$, that is, 
$$D_i=\cases{A_i&if $i\neq p$\cr
A_p-\sum_{m\in F^0}A_m+\sum_{m\in F^0}E_m&if $i=p.$\cr}
$$
The equation 
$$\eqalign{D_p&=A_p-\sum_{m\in F^0}A_m+\sum_{m\in F^0}E_m\cr
&=\sum_{m\in M}A_m+\sum_{k\in K}E_k-\sum_{m\in F^0}A_m+\sum_{m\in 
F^0}E_m\cr
&=\sum_{m\in H^0}D_m+\sum_{k\in K\cup F^0}E_k\cr}
$$
determines a primitive transfer of $D$ at $p$.  This primitive transfer 
is $B$ since if $i\neq p$ then $B_i=A_i=D_i$ and 
$$\eqalign{B_p&=A_p-\sum_{m\in M}A_m+\sum_{m\in M}E_m\cr
&=A_p-\sum_{m\in F^0}A_m-\sum_{m\in H^0}A_m+\sum_{m\in F^0}E_m+\sum_{
m\in H^0}E_m\cr
&=D_p-\sum_{m\in H^0}D_m+\sum_{m\in H^0}E_m.\cr}
$$
The statement about the graphs of the primitive transfers is obvious.  
\eop

\uj
We are now in position to show our main result.

\theorem If the matrix $B$ is a primitive transfer of $A$ then there is a 
sequence of matrices $A=C_1,\dots,C_i,\dots,C_n=B$ such that for all 
$1\le i\le n-1$ either $C_{i+1}$ is a size $1$ primitive transfer of $
C_i$ or $C_i$ is a 
size $1$ primitive transfer of $C_{i+1}$.  

\proof By Lemma \cite{complemma} there is a sequence of matrices 
$A=D_1,\dots,D_i,\dots,D_m=B$ such that for all $1\le i\le m-1$, $
D_{i+1}$ is the 
primitive transfer of $D_i$ and the graph of the primitive transfer is 
connected.  Inductively applying Lemma \cite{edgelemma} we can 
transform $D_i$ to $D_{i+1}$ for all $i$, using size $1$ primitive transfers.  By 
the note after Lemma \cite{edgelemma}, we never run out of non-loop 
edges.  
\eop\ 

\references
\uj
\i{[DS]}{D. Drinen and N. Sieben}{$C^{*}$-equivalences of graphs}{
J. Operator Theory (to appear)}       
\i{[EFW]}{M. Enomoto, M. Fujii and Y. Watatani}{$K_0$-groups and 
classifications of Cuntz-Krieger algebras}{Math. Japon. \bf 26 
\rm (1981), 443--460.}
\i{[Fra]}{J. Franks}{Flow equivalence of subshifts of finite type}{Ergo. 
Th. \& Dynam. Sys. \bf 4 \rm (1984), 53--66.}
\end